\title{On a sufficient condition that  $\sqrt{s}$ is simply normal to base 2, for $s$ not a perfect square}
\author{Richard Isaac}
\documentclass[11pt]{article}
\newtheorem{theorem}{Theorem}

\newtheorem{lemma}{Lemma}
\newtheorem{proposition}{Proposition}

\date{   }
\begin{document}
\maketitle
\begin{abstract}In \cite{ri:sn} the author introduced a condition, Condition (TU), and proved that its validity implies the simple normality to base 2 of $\sqrt{s}$, for $s$ not a perfect square. The argument also given in \cite{ri:sn} that Condition (TU) is indeed valid was cumbersome. We give here a simpler direct proof that Condition (TU) is true.    \footnote{
{\it \hspace*{1em} Keywords and phrases}. normal number, tail function, calculus of finite differences.}\\
\end{abstract}
 
\section{Introduction} 
In \cite{ri:sn} the author introduced a condition called Condition (TU) and proved that it implied the simple normality to base 2 of  $\sqrt{s}$ for $s$ not a perfect square. Also given was an argument that Condition (TU) is true. This argument was unnecessarily long, and was hard to follow according to some readers. Recently I have found a simpler proof of the validity of Condition (TU); it is presented in Theorem~\ref{final}. 

 Consider numbers $\omega$ in the unit interval, and represent the dyadic expansion of $\omega$ as
\begin{equation}
\omega=.x_{1}x_{2} \cdots ,\hspace{.4in} x_{i}=0 \mbox{ or } 1.  \label{expom}
\end{equation}
Also of interest is the dyadic expansion of $\nu=\omega^{2}$:
\begin{equation} \label{xu}
 \nu=\omega^{2}=.u_{1}u_{2} \cdots, \hspace{.4in} u_{i}=0 \mbox{ or } 1. \end{equation}
Throughout this paper it will be assumed that $\nu$ is irrational. Then $\omega$ is also irrational and both expansions are uniquely defined. It will be convenient to refer to the expansion of $\omega$ as an $x$ sequence and the expansion of $\nu$ as a $u$ sequence. A point of the unit interval can also be denoted by its coordinate representation, that is, $\omega=(x_{1}, x_{2}, \cdots)$ or $\nu=(u_{1}, u_{2}, \cdots)$.  The coordinate functions $X_{n}(\omega)=x_{n}$ and $U_{n}(\nu)=u_{n}$ give the $n$th coordinates of $\omega$ and $\nu$ respectively. 

Given any dyadic expansion $.s_{1}s_{2}\cdots$ and any positive integer $n$, the sequence of digits $s_{n}, s_{n+1}, \cdots$ is called a {\it tail} of the expansion. Two expansions are said to have the same tail if there exists $n$ so large that the tails of the sequences from the $n$th digit are equal.

The average
\begin{equation} \label{aver}
f_{n}(\omega)= \frac{x_{1}+x_{2}+\cdots+x_{n}}{n} \end{equation}
is the relative frequency of 1's in the first $n$ digits of the expansion of $\omega$. Simple normality for $\omega$ is the assertion that $f_{n}(\omega) \rightarrow 1/2$ as $n$ tends to infinity.
 Let $n_{k}$ be any fixed subsequence and define \begin{equation} \label{ff}
f(\omega)= \limsup_{k \rightarrow \infty} f_{n_{k}}(\omega). \end{equation}
We note that  
the function $f$ is a tail function with respect to the $x$ sequence, that is, $f(\omega)$ is determined by any tail $x_{n}, x_{n+1}, \cdots$ of the coordinates of \setcounter{footnote}{0}$\omega$.\footnote{ In fact, $f$ satisfies a more stringent requirement: it is an invariant function (with respect to the $x$ sequence) in the following sense: let $T$ be the 1-step shift transformation on $\Omega$ to itself given by
\[T(.x_{1}x_{2} \cdots) =.x_{2}x_{3} \cdots. \]
A function $g$ on $\Omega$ is invariant if $g(T\omega)=g(\omega)$ for all $\omega$. Any invariant function is a tail function.} 

We now observe that the average $f_{n}$, defined in terms of the $x$ sequence, can also be expressed as a function $h_{n}(\nu)$ of the $u$ sequence because the $x$ and $u$ sequences uniquely determine each other. This relationship has the simple form  
$f_{n}(\omega)=f_{n}(\sqrt{\nu})=h_{n}(\nu)$. Define $h(\nu)= \limsup_{k} h_{n_{k}}(\nu)$; then clearly 
$f(\omega)=h(\nu)$. \\
\noindent {\bf Definition}: Let $f$ be defined as in relation~\ref{ff} for any fixed subsequence $n_{k}$.  
We say that Condition (TU) is satisfied if $f(\omega)=h(\nu)$ is a tail function {\it with respect to the $u$ sequence} whatever the sequence $n_{k}$, that is, for any $\omega$ and any positive integer $n$,  $f(\omega)$ only depends on $u_{n}, u_{n+1}, \cdots$, the tail of the expansion of $\nu=\omega^{2}$.
(The notation ``TU'' is meant to suggest the phrase ``tail with respect to the $u$ sequence''.)

An immediate consequence of Condition (TU) is:
\begin{proposition} Let $\eta$ be the dyadic expansion of an irrational number. Let $\eta_{1}$ be a dyadic expansion that agrees with $\eta$ at all but a finite number of indices. If Condition (TU) is satisfied then 
\[\lim_{n} (f_{n}(\sqrt{\eta})-f_{n}(\sqrt{\eta_{1}}))=0. \]  \end{proposition}

\section{Proof of Condition (TU)} \begin{lemma} \label{one}
Let $\omega^{2}=\nu$.\\
(a) Let $u_{1}, u_{2}, \cdots, u_{r}$ be the initial segment of length $r$ of $\nu$. Then there exists a positive integer $N=N(\omega, r)$ such that each $u_{i},\, i \leq r$ is a function of $x_{1}, x_{2}, \cdots, x_{N}$. \\
(b)  Let $x_{1}, x_{2}, \cdots, x_{n}$ be the initial segment of length $n$ of $\omega$. Then there exists a positive integer $m=m(\nu,n)$ such that each $x_{j},\, j \leq n $ is a function of  $u_{1}, u_{2}, \cdots, u_{m}$.
\end{lemma}
The proof can be found in \cite{ri:sn}, lemmas 2 and 3.

The following arguments will use some elementary ideas from the calculus of finite differences. An introduction to these ideas may be found, for example, in \cite{sg:de}. We review some of the notation. Let $v(y_{1}, \cdots, y_{l})=v(\mbox{\boldmath $y$})$ be a function on the $l$-fold product space $S^{l}$ where the $y_{i} \in S$, a set of real numbers. 
Suppose that the variable $y_{i}$ is changed by the amount $\Delta y_{i}$ such that the $l$-tuple $\mbox{\boldmath $y^{(1)}$}=(y_{1}, \cdots, y_{l})$ is taken into $\mbox{\boldmath $y^{(2)}$}=(y_{1}+\Delta y_{1}, \cdots, y_{l}+\Delta y_{l})$ in the domain of definition of $v$. Put $v(\mbox{\boldmath $y^{(2)}$})-
v(\mbox{\boldmath $y^{(1)}$})= \Delta v$, and let 
\begin{eqnarray} \label{decomp1}
 \Delta v_{i}&=& v(y_{1}, \cdots, y_{i-1},y_{i}+\Delta y_{i}, y_{i+1}+\Delta y_{i+1, \cdots},
y_{l}+\Delta y_{l}) \hspace{.25in}\\
&-&v(y_{1}, \cdots, y_{i-1},y_{i}, y_{i+1}+\Delta y_{i+1, \cdots},
y_{l}+\Delta y_{l}). \nonumber \end{eqnarray} 
Then $\Delta v= \sum_{i} \Delta v_{i}$ is the total change in $v$ induced by changing all of the $y_{i}$, where this total change is written as a sum of step-by-step changes in the individual   
$y_{i}$. Formally, by dividing, we can write \begin{equation} \label{decomp2}
\Delta v= \sum_{i} (\Delta v_{i}/\Delta y_{i})\cdot \Delta y_{i}. \end{equation}
If some $\Delta y_{i_{0}}=0$, its coefficient in relation \ref{decomp2} has the form 0/0. 
No matter how the coefficient is defined in this case the contribution of the $i_{0}$ term to $\Delta v $ is 0. For our purposes it is convenient to define the coefficient to be $\Delta v_{i_{0}}$ evaluated as though $y_{i_{0}}$ were equal to 0 and $\Delta y_{i_{0}}$ were equal to 1.  

Let us then formally define the {\it partial difference of $v$ with respect to $y_{i}$, evaluated at the pair (\mbox{\boldmath $y^{(1)}$},\,\mbox{\boldmath $y^{(2)}$})} by \begin{eqnarray}  \label{denom}
\frac{\Delta v}{\Delta y_{i}} &=& \Delta v_{i}/ \Delta y_{i}, \hspace{1em} \mbox{ if } \Delta y_{i} \neq 0,\\
& = & \Delta v_{i} \mbox{ evaluated as though }  y_{i}=0 \mbox{ and } \Delta y_{i}=1,\hspace{1em} \mbox{ if } \Delta y_{i}=0 .\nonumber \end{eqnarray}
Notice that the forward slash (/) in this relation expresses division and the horizontal slash on the left hand side is the partial difference operator. 

The sum $\Delta v$ of relation~\ref{decomp2} is called the {\it total difference of $v$} evaluated at the given pair and can now be written
\begin{equation} \label{decomp3}
\Delta v= \sum_{i} \frac{\Delta v}{\Delta y_{i}}\cdot \Delta y_{i}. \end{equation}
The $i$th summand in relation~\ref{decomp3} is called the {\it $i$th partial difference} of $v$ relative to the given pair.
The partial and total differences are the discrete analogs of  the partial  and total differentials in the theory of differentiable functions of several real variables and the partial difference with respect to a given $y$ variable is the analog of the partial derivative. The $i$th partial difference of $v$  at a given pair is a measure of the contribution of $\Delta y_{i}$ to $\Delta v$ when all the other $y$ variables are held constant.   

Returning to our particular problem, we say that $\omega$ and $\nu=\omega^{2}$ are points (or expansions) that {\it correspond} to one another. As seen in Section 1 the average $f_{n}(\omega)$ of relation~\ref{aver} can be written as a function $h_{n}(\nu)$.  With a slight abuse of notation we can write
 \begin{equation} \label{xy2}
f_{n}(x_{1}, \cdots, x_{n})= f_{n}(\omega )=h_{n}(\nu)=h_{n}(u_{1}, u_{2}, \cdots). \end{equation}
Fix a point $\omega$ with corresponding point $\nu$,  and for each $x_{j}$ let $\Delta x_{j}$ be a given increment chosen independently ($\Delta x_{j}=0, 1, $ or $-1$). Let $\omega^{(1)}$ have coordinates $x_{j}+\Delta x_{j}$ and let $\nu^{(1)}$ correspond to $\omega^{(1)}$. Let the $i$th coordinate of  $\nu^{(1)}$ be $u_{i}+\Delta u_{i}$.  Thus the changes $\Delta x_{j}$ in the $x$ coordinates have induced changes $\Delta u_{i}$ in the $u$ coordinates. Of course this process could have been reversed: independent changes in the $u$'s  induce changes in the $x$'s.

The following two lemmas are  finite difference analogs of the total differential formulas in the theory of differentiable functions of a function of several variables.  The first result is fairly evident.  
\begin{lemma} \label{ecks}  At the pair  $(\omega, \omega^{(1)}),\,  \Delta f_{n}$ can be represented as a total difference 
\begin{eqnarray} \label{fx} 
\Delta f_{n}&=&f_{n}( \omega^{(1)})-f_{n}(\omega)=\\
&&  f_{n}(x_{1}+\Delta x_{1}, x_{2}+\Delta x_{2}, \cdots x_{n}+\Delta x_{n}) - f_{n}(x_{1}, x_{2}, \cdots x_{n}) \nonumber \\
&& =\frac{1}{n} \sum_{ 1 \leq j \leq n} \Delta x_{j} \nonumber \end{eqnarray} \end{lemma}
Proof: Decompose according to the recipe given in relations~\ref{decomp1} to \ref{decomp3} to \setcounter{footnote}{0}get \footnote{Our definitions  require the ``denominator'' of a partial difference to be a variable, so strictly speaking $x_{j}$ in this relation should be replaced by $X_{j}$, the $j$th coordinate variable, with an added notation that it is evaluated at the given base point $\omega$. The present notation is simpler and will be followed throughout.}
 \[\frac{\Delta f_{n}}{\Delta  x_{j}}=\frac{1}{n}, \hspace{1em} j \leq n \mbox{ and } =0, \hspace{1em} j>n.\]
The next lemma is more interesting.     
\begin{lemma} \label{fini}   At the pair $(\nu, \nu^{(1)})$, $\Delta h_{n}$ can be represented as a total difference
\begin{eqnarray} \label{delh}
\Delta h_{n}&=&h_{n}(\nu^{(1)})- h_{n}(\nu) = h_{n}(u_{1}+\Delta u_{1}, u_{2}+\Delta u_{2}, \cdots)-h_{n}(u_{1}, u_{2}, \cdots)\nonumber \\& & = \sum_{i \geq 1} \frac{\Delta h_{n}}{\Delta u_{i}}\, \Delta u_{i}=          \sum_{i \geq 1} (\Delta h_{n,i}/\Delta u_{i})\, \Delta u_{i}, \hspace{1em}  \Delta u_{i} \neq 0\end{eqnarray}  
where
\begin{eqnarray} \label{pd}  \Delta h_{n,i}=\Delta h_{n,i}(\nu, \nu^{(1)})=&&\\
h_{n}(u_{1}, \cdots, u_{i-1}, u_{i}+\Delta u_{i}, u_{i+1}+\Delta u_{i+1}, \dots)-&& \nonumber\\
h_{n}(u_{1}, \cdots, u_{i-1}, u_{i}, u_{i+1}+\Delta u_{i+1}, \dots)&&. \nonumber \end{eqnarray}
The formally infinite sum of relation~\ref{delh} reduces to a finite sum. More precisely, given the pair $(\nu, \nu^{(1)})$, there exists an integer $m$ such that the partial differences  $\Delta h_{n,i}/\Delta u_{i}=0$ for all $i>m$. The number of non-vanishing terms in the sum depends on $\nu$ and  $n$.
\end{lemma}
Proof: The recipe given in relation~\ref{delh} for decomposing $\Delta h_{n}$ is given by the definitions stated in relations~\ref{decomp1} through \ref{decomp3}.  The pair $(\nu, \nu^{(1)})$ corresponds to the pair $(\omega, \omega^{(1)})$. To see that the sum in relation~\ref{delh} is finite, note that the function $h_{n}=f_{n}$ only depends on $x_{1}, x_{2}, \cdots x_{n}$. Given $\nu$, Lemma~\ref{one} proves the existence of an integer $m$ such that for all $i>m$
\[ \frac{\Delta x_{j}}{ \Delta u_{i}}=0, \hspace{1.5em} 1 \leq j \leq n. \]
 and therefore the terms $\Delta h_{n,i}$ of relation~\ref{pd} are 0 for $i>m$. Thus the terms in the sum of relation~\ref{delh} vanish for $i>m$ and the formula of relation~\ref{delh} represents a finite sum. This concludes the proof of the lemma. 

We have seen (lemma~\ref{ecks}) that the partial difference of $f_{n}$ with respect to any fixed $x_{j}$ is $1/n $. This means that the contribution to changes in the averages $f_{n}$ of any change in a single $x_{j}$  tends to 0. But how about the partial differences of $h_{n}$ with respect to a single fixed $u_{i}$? Does a change in $u_{i}$ induce changes in $h_{n}$ that die out in the limit? The answer is not obvious. Let $u_{r}$ be a fixed $u$ variable. Relation~\ref{pd} shows that \begin{equation} \label{llte}
\Delta h_{n,r}= \frac{1}{n} \sum_{ 1 \leq j \leq n} \Delta x_{j}^{\prime}  \end{equation}
for some sequence of changes of $x$ variables (see proof of theorem~\ref{final} below). If $\Delta u_{r} \neq 0$ the partial difference of $h_{n}$ with respect to $u_{r}$ just differs from $\Delta h_{n,r}$ by a factor of $\pm 1$. Therefore the questions posed above reduce to asking what the limit points of the right hand side of relation~\ref{llte} are. Heuristic considerations suggest why we might expect the averages in relation~\ref{llte} to converge to 0:  $h_{n}$ depends on larger and larger initial segments of $u$ variables as $n$ increases and a certain symmetry exists in the problem. It seems reasonable to suspect that change in a single $u$ variable is not going to have much of an effect on $h_{n}$ for large $n$. We now set out to prove that this suspicion is true. 

Since the $u$ variables are functions of the $x$ variables and vice versa, it is possible to consider either  set of variables independent and the other set dependent on them. We choose to take the $x$ variables independent. The power of this approach becomes apparent in the next result which solves our problem.
 \begin{theorem}  \label{final} Assume that the $u$ variables are functions of independent $x$ variables. Then \\
 (a):  For all $i$, the partial differences of $h_{n}$ with respect to $u_{i}$  in relation~\ref{delh} satisfy 
 \begin{equation} \label{limm}
\lim_{n} \frac{\Delta h_{n}}{\Delta u_{i}}=0. \end{equation} 
(b): Condition (TU) is true.  \end{theorem}
Proof of (a):  The $i$th partial differences $\Delta u_{i}$ referenced in (a) are all non-zero, so let $r$ be a fixed positive integer with $\Delta u_{r} \neq 0$. Consider relation~\ref{pd}. The right hand side expresses $\Delta h_{n, r}$ as the difference  $h_{n}(\nu_{2})- h_{n}(\nu_{1})$ evaluated at the two points 
\begin{eqnarray} \label{nu12}
\nu_{2}&=&(u_{1}, \cdots, u_{r-1}, u_{r}+\Delta u_{r}, u_{r+1}+\Delta u_{r+1}, \dots)  \mbox{ and }  \\
\nu_{1}&=&(u_{1}, \cdots, u_{r-1}, u_{r}, u_{r+1}+\Delta u_{r+1}, \dots)    
 \nonumber  \end{eqnarray}
The irrationality of $\nu^{(1)}$ implies that $\nu_{1}$ and $\nu_{2}$ are also irrational.
For $k=1,2, $ let $\nu_{k}$ correspond to $\omega_{k}=(x_{(1,k)}, x_{(2,k)}, \cdots)$ and put 
$x_{(j,2)}-x_{(j,1)}= \Delta  x^{\prime}_{j}$. Let the differences in the $u$ coordinates at $(\nu_{1}, \nu_{2})$ be denoted by $\Delta u_{i}^{\prime}$. Then $\Delta u_{i}^{\prime}=0$ for $i \neq r,\, \Delta u_{r}^{\prime}=\Delta u_{r}$. We study the functions $f_{n}$ and $h_{n}$ at the pairs $(\omega_{1},\omega_{2})$ and $(\nu_{1},\nu_{2})$ respectively.    At the pairs $(\omega_{1},\omega_{2})$  and  $(\nu_{1},\nu_{2})$,  lemmas~\ref{ecks} and   \ref{fini} correspond  to    \begin{equation} \label{hf}
   \Delta h_{n}^{\prime}=  \Delta f_{n}^{\prime}=f_{n}(\omega_{2})-f_{n}(\omega_{1})= \frac{1}{n} \sum_{ 1 \leq j \leq n} \Delta x_{j}^{\prime}  \end{equation}
and  \begin{equation} \label{vers1}
\Delta h_{n}^{\prime}= h_{n}(\nu_{2})-h_{n}(\nu_{1})=\frac{\Delta h_{n}^{\prime}}{\Delta u_{r}}\, \Delta u_{r},    \mbox{ where }   \hspace{1em}     \frac{ \Delta h_{n}^{\prime}}{\Delta u_{r}}= \pm \frac{1}{n} \sum_{ 1 \leq j \leq n} \Delta x_{j}^{\prime}.       \end{equation}
At the pair $(\nu_{1},\nu_{2})$ 
\[\Delta h_{n}^{\prime}= \frac{\Delta h_{n}^{\prime}}{\Delta u_{r}} \Delta u_{r}=h_{n}(\nu_{2})-h_{n}( \nu_{1}) =\frac{\Delta h_{n}}{\Delta u_{r}} \Delta u_{r}=\Delta h_{n, r}  \]
so that \begin{equation} \label{ssame}
\frac{\Delta h_{n}^{\prime}}{\Delta u_{r}}=\frac{\Delta h_{n}}{\Delta u_{r}} \end{equation}
and  \begin{equation} \label{vers2}
\Delta h_{n}^{\prime}= \frac{\Delta h_{n}}{\Delta u_{r}}\, \Delta u_{r}. \end{equation}
Since $u_{r}$ is a function of the $x$ variables, at the pair $(\omega_{1}, \omega_{2})$  relations~\ref{decomp1} through \ref{decomp3} give the representation
\[\Delta u_{r}= \sum_{j \geq 1} \frac{\Delta u_{r}}{\Delta x_{j}^{\prime}}\, \Delta x_{j}^{\prime}. \]
By lemma~\ref{one} there exists $N=N(\omega_{1},r)$ such that the changes $ \Delta x_{j}^{\prime},\, j>N$ cause no change in $u_{r}$, that is,
\[ \frac{\Delta u_{r}}{\Delta x_{j}^{\prime}}=0,  \hspace{1.5em} j>N. \]
 It follows that there is the finite decomposition \begin{equation} \label{uxx}
 \Delta u_{r}= \sum_{1 \leq j \leq N} \frac{\Delta u_{r}}{\Delta x_{j}^{\prime}}\, \Delta x_{j}^{\prime}. \end{equation}
Using relation~\ref{uxx}, relation~\ref{vers2} can be rewritten    \begin{equation} \label{vers3}
   \Delta h_{n}^{\prime} =    \sum_{1 \leq j \leq N} \frac{\Delta h_{n}}{\Delta u_{r}} \frac{\Delta u_{r}}{\Delta x_{j}^{\prime}}\, \Delta x_{j}^{\prime}.  \end{equation}
Let $n_{k}$ be any subsequence for which  there is convergence in relation~\ref{vers3}, that is, 
\[ \lim_{k} \Delta h_{n_{k}}^{\prime}=  \Delta h^{\prime} \mbox{ and } \lim_{k} \frac{\Delta h_{n_{k}}}{\Delta u_{r}}=  \frac{\Delta h}{\Delta u_{r}} \]
where the right hand sides are defined by the existing limits. Then
 \begin{equation} \label{vers4}
   \Delta h^{\prime} =    \sum_{1 \leq j \leq N} \frac{\Delta h}{\Delta u_{r}} \frac{\Delta u_{r}}{\Delta x_{j}^{\prime}}\, \Delta x_{j}^{\prime}.  \end{equation}
Let $p \leq N$ be an index with \begin{equation} \label{rp}
\frac{\Delta u_{r}}{ \Delta x_{p}^{\prime}}\,\Delta x_{p}^{\prime} \neq 0. \end{equation}
Such $p$ exists by relation~\ref{uxx} since $\Delta u_{r} \neq 0$.
Now observe that  $\Delta h^{\prime}=  \lim_{k} \Delta h^{\prime}_{n_{k}}$  is a tail function considered as a function of the $ \Delta x_{j}^{\prime}$ variables (see relation~\ref{hf}),  so is not a function of  $\Delta x_{j}^{\prime}$  for any fixed $j$.  Moreover, \begin{equation} \label{nor}
 \frac{\Delta h}{\Delta u_{r}}  \mbox{ is a tail function with respect to the $\Delta x_{j}^{\prime}$  (see relations~\ref{vers1} and \ref{ssame})}. \end{equation}
Take the partial difference with respect to the $p$th coordinate variable on both sides of relation~\ref{vers4}. By the tail property of  $\Delta h^{\prime}$ stated above,  
\begin{equation} \label{jp0}
\frac{\Delta h^{\prime}}{\Delta x^{\prime}_{p}}=0. \end{equation}
Independence of the $x$ variables implies
 \[ \frac{\Delta x_{j}^{\prime}}{\Delta x_{p}^{\prime}}=0, \hspace{2em} j \neq p.    \]
 Use relations~\ref{rp} and \ref{nor} and the foregoing relation to see that the partial difference with respect to the $p$th coordinate variable on the right hand side of relation~\ref{vers4} can be written
 \begin{equation} \label{vers5}
      \sum_{1 \leq j \leq N} \frac{\Delta h}{\Delta u_{r}} \frac{\Delta u_{r}}{\Delta x_{j}^{\prime}}\, \frac{\Delta x_{j}^{\prime}}{\Delta x_{p}^{\prime}}= \frac{\Delta h}{\Delta u_{r}} \frac{\Delta u_{r}}{\Delta x_{p}^{\prime}}= \pm \frac{\Delta h}{\Delta u_{r}}  \end{equation}
      and so  from relations~\ref{jp0} and \ref{vers4} relation~\ref{vers5} implies \begin{equation} \label{imp}
      \frac{\Delta h}{\Delta u_{r}}=0. \end{equation}
 The subsequence $n_{k}$ is associated with an arbitrary limit point so the above argument shows this limit point is unique, that is 
\[  \lim_{n}  \frac{\Delta h_{n}}{\Delta u_{r}}=0 \]
and this proves (a).\\
Proof of (b): Given any subsequence $n_{k}$ and any positive integer M, part (a) of this theorem proves that in relation~\ref{delh}
\[\limsup_{k} \Delta h_{n_{k}}= \limsup_{k} \left(\sum_{i >M} \frac{\Delta h_{n_{k}}}{\Delta u_{i}}\, \Delta u_{i} \right).\]
The relation shows that $\limsup_{k} \Delta h_{n_{k}}$ does not depend on the differences of any initial segment of $u$ coordinates for the given pair in relation~\ref{delh}. Since lemma~\ref{fini} makes no restrictions on pairs (other than they are well defined), this assertion is true for all meaningful pairs. This implies that  $\limsup_{k}  h_{n_{k}}=\limsup_{k} f_{n_{k}}$ is a tail function with respect to the $u$ variables, that is, Condition (TU) is true.  

\vspace{.4in}

{\bf Emeritus Professor}  

{\bf Lehman College and Graduate Center, CUNY}

{\bf {\it email:} richard.isaac@lehman.cuny.edu}
\end{document}